\documentclass{amsart}

\usepackage{amsmath} 
\usepackage{amssymb}
\usepackage{mathrsfs}

\newtheorem{theorem}{Theorem}

\theoremstyle{definition}

\theoremstyle{remark}
\newtheorem{remark}[theorem]{Remark}

\newcommand{\rest}{{\restriction}}

\newcommand{\mn}{{\medskip\noindent}}
\newcommand{\sn}{{\smallskip\noindent}}

\newcommand{\bbR}{{\mathbb R}}

\newcount\skewfactor
\def\mathunderaccent#1#2 {\let\theaccent#1\skewfactor#2
\mathpalette\putaccentunder}
\def\putaccentunder#1#2{\oalign{$#1#2$\crcr\hidewidth
\vbox to.2ex{\hbox{$#1\skew\skewfactor\theaccent{}$}\vss}\hidewidth}}

\newenvironment{PROOF}[2][\proofname.]
   {\begin{proof}[#1]\renewcommand{\qedsymbol}{\textsquare$_{\rm #2}$}}
   {\end{proof}}

\begin{document}

\title {Inner product space with no ortho-normal basis without choice}
\author {Saharon Shelah}
\address{Einstein Institute of Mathematics\\
Edmond J. Safra Campus, Givat Ram\\
The Hebrew University of Jerusalem\\
Jerusalem, 91904, Israel\\
 and \\
 Department of Mathematics\\
 Hill Center - Busch Campus \\ 
 Rutgers, The State University of New Jersey \\
 110 Frelinghuysen Road \\
 Piscataway, NJ 08854-8019 USA}
\email{shelah@math.huji.ac.il}
\urladdr{http://shelah.logic.at}
\thanks{The author thanks Alice Leonhardt for the beautiful typing.\\ 
   Publication E68.}

\date{August 3, 2010}

\begin{abstract}
We prove in ZF that there is an inner product space, in fact, nicely
definable with no orthonormal basis.
\end{abstract}

\maketitle
\numberwithin{equation}{section}

The theorem below is known in ZFC, but probably not in ZF; really we
use the simple black box (see \cite{Sh:309}).
\begin{theorem}  
\label{a1}
\underline{(ZF)}
There is an inner-product space $V$ over $\bbR$ with no ortho-normal basis.
\end{theorem}

\begin{remark}
In fact, nicely definable one, here - Borel.
\end{remark}

\begin{PROOF}{\ref{a1}}  
\smallskip

\noindent
\underline{Stage A}:

Let $V_1$ be the Hilbert space over $\bbR$ with orthonormal basis
$\{x_\eta:\eta \in {}^{\omega \ge}\omega\}$, so an element $x$ has a
unique representation as $x = \Sigma\{a_{x,\eta} x_\eta:\eta \in
{}^{\omega >}\omega\}$ with $a_{x,\eta} \in \bbR$ and norm $< \infty$
so supp$_1(x) := \{\eta:a_{x,\eta} \ne 0\}$ is countable and
supp$^1_k(x) := \{\eta:|a_\eta| \ge \frac{1}{k+1}\}$ finite for every $k
< \omega$ where the norm is $\Sigma\{a^2_{x,\eta}:\eta \in
{}^{\omega >}\omega\}$.  The inner product is $((\Sigma a_\eta
x_\eta),(\Sigma a'_\eta x_\eta)) = \Sigma\{a_\eta a'_\eta:\eta \in
{}^{\omega \ge}\omega\} \in \bbR$.

For $\eta \in {}^{\omega}\omega$ let $y_\eta = x_\eta + \sum\limits_{n
< \omega} \frac{1}{2^n} x_{\eta \rest n}$.

Let $V$ be the subspace of $V_1$ generated by $\{x_\eta:\eta \in
{}^{\omega >}\omega\} \cup \{y_\eta:\eta \in {}^{\omega}\omega\}$ so
as a vector space it is $\bigoplus\limits_{\eta \in {}^{\omega
>}\omega} \bbR x_\eta \oplus \bigoplus\limits_{\eta \in
{}^{\omega}\omega} \bbR y_\eta$ and it ``inherits" the innder product
from $V_1$.

Toward contradiction assume that $\{z_s:s \in S\}$ is an ortho-normal
basis of $V$.  So every $x \in V$ has the unique representation 
$\sum\limits_{s \in S} b_{x,s} z_s$, where $b_{s,n} \in \bbR$ and 
for $k \in [1,\omega)$
and $x \in V$ let supp$^2_k(x) := \{s \in S:|b_{x,s}| \ge
\frac{1}{2k+1}\}$, so finite and supp$_2(x) := \{s \in S:b_{x,s} \ne
0\}$ so countable.
\medskip

\noindent 
\underline{Stage B}:

We choose $\eta_n$ by induction on $n$ such that:
\mn
\begin{enumerate}
\item[$\boxplus_1$]  $(a) \quad \eta_\eta \in {}^n \omega$
\sn
\item[${{}}$]   $(b) \quad \eta_m = \eta_n \rest m$ if $m < n$
\sn
\item[${{}}$]   $(c) \quad$ if $n = m+1$ then $\eta_n = \eta_m \char
94 \langle i \rangle$ with $i < \omega$ minimal such that: 

\hskip25pt if $\ell
\le m,s \in \text{ supp}^2_m(x_{\eta_\ell})$ and $\nu \in \text{
supp}^1_m(z_s) \subseteq {}^{\omega >}\omega$ then 

\hskip25pt $\neg(\eta_m
\char 94 \langle i \rangle \trianglelefteq \nu)$.
\end{enumerate}
\mn
This is well defined as in clause (c), supp$^2_m(x_{\eta_\ell})$ is a
finite subset of $S$ and for each $s \in \text{
supp}^2_m(x_{\eta_\ell})$, the set supp$^1_m(z_s)$ is a finite subset
of ${}^{\omega >}\omega$.

Lastly, let
\mn
\begin{enumerate}
\item[$\boxplus_2$]  $(a) \quad \eta_\omega := \cup\{\eta_n:n < \omega\} \in
{}^\omega \omega$
\sn
\item[${{}}$]   $(b) \quad S_1 = \cup\{\text{supp}_2(x_\rho):\rho
\triangleleft \eta_\omega\}$
\sn
\item[${{}}$]   $(c) \quad S_2 = S \backslash S_1$
\sn
\item[${{}}$]   $(d) \quad X_\ell$ is the closure inside $V$ of
$\oplus\{\bbR z_s:s \in S_\ell\}$ for $\ell=1,2$
\sn
\item[${{}}$]   $(e) \quad S_{1,n} :=
\cup\{\text{sup}^2_m(x_{\eta_\ell}):m,\ell \le n\}$.
\end{enumerate}
\mn
Note
\mn
\begin{enumerate}
\item[$\boxplus_3$]  $V = X_1 \oplus X_2$, i.e. $X_1,X_2$
are orthogonal but $X_1 + X_2$ is $V$
\sn
\item[$\boxplus_4$]  $S_1 = \cup\{\text{supp}^2_m(x_{\eta_n}):n <
\omega,m < n\} = \cup\{S_{1,n}:n < \omega\}$
\sn
\item[$\boxplus_5$]   $\eta_n \in S_1$ for $n < \omega$.
\end{enumerate}
\medskip

\noindent
\underline{Stage C}:  As $y_{\eta_\omega} \in V$ see Stage A and
$\boxplus_2(a)$ of Stage B, recalling $\boxplus_3$
\mn
\begin{enumerate}
\item[$\otimes_1$]   there are $y^1 \in X_1,y^2 \in X_2$ such that
$y_{\eta_\omega} = y^1 + y^2$.
\end{enumerate}
\mn
Also
\mn
\begin{enumerate}
\item[$\otimes_2$]   $\{\rho:\eta_{n+1} \trianglelefteq 
\rho \in {}^{\omega \ge}
\omega\}$ is disjoint to $\cup\{\text{supp}^1_m(z_s):s \in S_{1,n}\}$
for every $n < \omega$.
\end{enumerate}
\mn
[Why? By the choice of $\eta_{n+1}$ in $\boxplus_1(c)$.]
\mn
\begin{enumerate}
\item[$\otimes_3$]  $\eta_\omega \notin \text{ supp}_1(z_s) =
\cup\{\text{supp}^1_m(z_s):m < omega\}$ for every $s \in S_1$.
\end{enumerate}
\mn
[Why?  The $\notin$ by $\otimes_2$.]

Hence by $\boxplus_3$
\mn
\begin{enumerate}
\item[$\otimes_4$]   if $s \in S_1$ then $y_{\eta_\omega},z_s$ are orthogonal
(in $V_1$).
\end{enumerate}
\mn
But
\mn
\begin{enumerate}
\item[$\otimes_5$]  $(y_{\eta_\omega},x_{\eta_n}) = \frac{1}{2^n}$.
\end{enumerate}
\mn
[Why?  By the choice of $y_\eta$ is stage N.]  

By $\boxplus_5 +
\otimes_4 + \otimes_5$ we get contradiction.
\end{PROOF}


\begin{thebibliography}{}

\bibitem[Sh:309]{Sh:309}
Saharon Shelah, \emph{{Black Boxes}}, {}, 0812.0656. 0812.0656. 0812.0656.

\end{thebibliography}
\end{document}